\documentclass[11pt,reqno]{amsart}
\usepackage[utf8]{inputenc}
\usepackage{setspace}
\usepackage{geometry}
\usepackage{enumerate}
\usepackage{enumitem, xcolor, amssymb,latexsym,amsmath,bbm}
\usepackage{mathtools}
\usepackage[mathscr]{euscript}
\usepackage{amsmath}
\usepackage{amssymb}
\usepackage{enumitem} 
\newcommand{\R}{\mathbb{R}}

\usepackage{tikz}
\usepackage{wrapfig}
\usepackage{enumerate}
\usepackage{graphicx}
\usepackage{subfigure}

\usetikzlibrary{arrows.meta}

\usetikzlibrary{decorations.markings}

\usepackage[colorlinks=true,citecolor=blue, linkcolor=blue,urlcolor=blue]{hyperref}

\calclayout

\theoremstyle{plain}
\newtheorem{theorem}{Theorem}[section]
\newtheorem{corollary}[theorem]{Corollary}
\newtheorem{lemma}[theorem]{Lemma}
\newtheorem{proposition}[theorem]{Proposition}
\theoremstyle{definition}
\newtheorem{definition}[theorem]{Definition}

\theoremstyle{remark}
\newtheorem{remark}[theorem]{Remark}

\numberwithin{equation}{section}
\numberwithin{figure}{section}
\allowdisplaybreaks

\title{Length of filling pairs on punctured surfaces}
\author{Bhola Nath Saha and Bidyut Sanki}

\date{\today}

\begin{document}


\begin{abstract}
    A pair $(\alpha, \beta)$ of simple closed curves on an oriented surface $S_{g,n}$ of genus $g$ and with $n$ punctures is called a filling pair if the complement of the union of the curves is a disjoint union of topological disks with at most one puncture. In this article, we study the lengths of filling pairs on once-punctured hyperbolic surfaces. In particular, we find a lower bound of the lengths of filling pairs. Furthermore, we show that this lower bound is the best one which depends only on the topology of the surface.
\end{abstract}

\maketitle

\tikzset{->-/.style={decoration={
  markings,
  mark=at position #1 with {\arrow{>}}},postaction={decorate}}}
  \tikzset{-<-/.style={decoration={
  markings,
  mark=at position #1 with {\arrow{<}}},postaction={decorate}}}

\section{Introduction}
Let $S_g$ denote a closed and orientable surface of genus $g$. A collection $\{ \alpha_i\mid i=1,\cdots,m\}$ of simple closed curves on $S_{g}$ is called a filling system of size $m$ if $S_{g}\setminus \cup_{i=1}^m\alpha_i$ is a disjoint union of topological disks. It is assumed that $\alpha_i$'s are pairwise homotopically distinct and in minimal position. The study of filling systems arises in the work of Thurston~\cite{thurstonspine} on constructing the spine of Teichm\"uller space of closed and oriented hyperbolic surfaces, where he has proposed that the subset $\chi_ g$ of the Teichm\"uller space consisting of those hyperbolic surfaces that admit a filling set of systolic geodesics as a candidate spine, i.e., there exists a mapping class group equivariant deformation retraction of the Teichm\"uller space onto $\chi_g$. Later, in \cite{Penner}, Penner used filing systems of surfaces to construct pseudo-Anosov homeomorphisms.

The process of filling systems on closed surfaces has undergone extensive research. Anderson--Parlier--Pettet~\cite{Anderson} have studied the size of a $k$-filling system on a closed surface, which is a filling system such that the simple closed curves pairwise intersect at most $k$-times. The authors have shown that the asymptotic growth rate for the minimal cardinality of such a filling system is $2 \sqrt{g}/\sqrt{k}$ as $g\to \infty$. They have also shown that the cardinality of a filling system of systoles in a closed hyperbolic surface of genus $g$ is at least $g/\ln{g}$.

Aougab--Huang \cite{Aougab} have studied minimally intersecting filling pairs, i.e., a filling system with two curves having minimal possible geometric intersection number. The mapping class group (Chapter 2 \cite{farb2011primer}) acts on the set of minimally intersecting filling pairs. Aougab--Huang have constructed exponentially many mapping class group orbits of such filling pairs. They have defined a function $\mathcal{F}_g$ on the moduli space which, given a hyperbolic metric $\sigma$, outputs the length of the shortest minimally intersecting filling pair for the metric $\sigma$ and shown that $\mathcal{F}_g$ is a topological Morse function. Furthermore, they have provided a sharp lower bound of $\mathcal{F}_g$, which depends only on the topology of the surface, and conjectured that the same lower bound would work if we consider any filling pairs instead of minimally intersecting.

Sanki--Vadnere have proved a generalized isoperimetric inequality for multiple hyperbolic polygons in Theorem 1.1~\cite{Sanki} and shown that the Aougab-Huang conjecture follows as a corollary of this isoperimetric inequality. Later, Gaster~\cite{Gaster} has given a direct and short proof of the Aougab--Huang conjecture.

For a surface $S_{g,n}$ of genus $g$ and $n$ punctures, a collection of simple closed curves is called a filling system if the complement of the union of the curves is a disjoint union of topological disks and once-punctured disks. Fanoni--Parlier \cite{Fanoni} have extended the study of filling systems to punctured surfaces. They have studied the cardinality of a $k$-filling system on the surface $S_{g,n}$. Furthermore, they have shown that a filling set of simple closed curves that intersect pairwise at most once contains at least $\sqrt{3n}$ curves. The authors have also shown that in a hyperbolic surface of signature $(g,n)$ and with systolic length $l$, a filling set of systolic geodesics has cardinality at least $[2\pi (2g-1)+\pi (n-2)]/4l$.

Motivated by the works of Aougab--Huang~\cite{Aougab}, Sanki--Vadnere~\cite{Sanki} and Gaster~\cite{Gaster}, in this article, we study the length function of filling pairs on once-punctured hyperbolic surfaces.

Let $(\alpha, \beta)$ be a filling pair on a surface $S_{g,n}$ of genus $g$ and with $n$ punctures and $\mathcal{T}_{g,n}$ be the Teichmuller space of $S_{g,n}$. The length of a filling pair $(\alpha, \beta)$ on a hyperbolic surface $X\in \mathcal{T}_{g,n}$ is defined by 
$$l_X(\alpha, \beta)=l_X(\alpha)+l_X(\beta),$$
where $l_X(\alpha)$ is the length of the unique simple closed geodesic in the free homotopy class of $\alpha$ on $X$ and similarly $l_X(\beta)$. Now, we define a function $\mathcal{Y}_{g,n}:\mathcal{T}_{g,n}\to \R$ by 
$$\mathcal{Y}_{g,n}(X)=\inf\{ l_X(\alpha,\beta)\mid (\alpha,\beta) \text{ is a filling pair on } X\}.$$
It is straightforward to see that $\mathcal{Y}_{g,n}$ is proper and a topological Morse function by the similar arguments as Aougab--Huang used in the proof of Theorem 1.3 \cite{Aougab}. Here, we find a global lower bound of the function $\mathcal{Y}_{g,1}$, which also depends only on the topology of the surface. In particular, we prove the following theorem.
\begin{theorem}\label{thm 1.1}
    The function $\mathcal{Y}_{g,1}$ satisfies 
    $$\mathcal{Y}_{g,1}(X)\geq (8g-4)\ln(\sqrt{2}+1), \text{ for all } X\in \mathcal{T}_{g,1}.$$
    Furthermore, this inequality is sharp.
\end{theorem}
The main ingredient of the proof of Theorem \ref{thm 1.1} is the isoperimetric inequality for the so-called polygonal cusp. A \emph{polygonal cusp} is a hyperbolic object that is topologically a once-punctured disk with a piecewise geodesic boundary of finite length (see Section 2 for more details). The geodesic segments are called \emph{sides} of the cusp. A Polygonal cusp is called \emph{regular} if all the side lengths are equal and the interior angles between two consecutive sides are also equal (for details, see Section 2). We prove the following version of the isoperimetric inequality. 

\begin{theorem}\label{lem: isom}
    Among all $p$-sided polygonal cusps with a fixed perimeter $l$, the regular one with interior angle 
    $$\theta = 2\sin^{-1}(\operatorname{sech}(l/2p)),$$
    has the maximum area.
\end{theorem}
The proof of Theorem~\ref{lem: isom} is based on basic hyperbolic trigonometry.
Furthermore, we note that the dual statement of Theorem~\ref{lem: isom} is also true. In particular, we prove that among all $p$-sided polygonal cusps with a fixed area $a$, the regular one has the least perimeter in Theorem~\ref{dual isom}. This dual statement is used in the proof of Theorem~\ref{thm 1.1}.
\subsection*{Structure of the paper}
In Section 2, we define polygonal cusps in detail, prove several technical lemmas, and establish the isoperimetric inequality stated in Theorem \ref{lem: isom}. In Section 3, we prove Theorem \ref{thm 1.1}.

\section{Isoperimetric inequality for polygonal cusp}
In this section, our goal is to prove Theorem~\ref{lem: isom}. We begin with defining polygonal cusps. Next, we develop two technical lemmas which are essential for the proof of the theorem. The lemmas deal with a version of the isoperimetric inequality for hyperbolic triangles. Then, we provide a proof of Theorem~\ref{lem: isom}. We conclude the section by proving the dual statement of the theorem in Theorem~\ref{dual isom}.

\subsection*{Polygonal cusp} Recall that a domain $D$ in a  non-compact hyperbolic surface is called a cusp if it is isometric to 
    $$X=\Tilde{X}/\Gamma,$$
    where $\Tilde{X}=\{ x+iy\in \mathbb{H} \mid y\geq 1\}$ and $\Gamma= \langle T\mid T(z)=z+1 \text{ for all } z\in \mathbb{H} \rangle < \mathrm{Isom^+(\mathbb{H})}$ (for details, see Example 1.6.8. \cite{Buser}).

Given a filling pair $(\alpha,\beta)$, where $\alpha$ and $\beta$ are simple closed geodesics on a punctured surface $X\in \mathcal{T}_{g,n}$, a component of $X\setminus (\alpha\cup\beta)$ is either a hyperbolic polygon or a polygon with a puncture. Motivated by this decomposition, we introduce the following construction.

\textbf{A construction:} Consider a collection of $m$ hyperbolic triangles, each having exactly one ideal vertex (for $m=5$, see the shaded triangles in Figure~\ref{fig: Cusp}(i)). Identify the sides of these triangles in pairs so that the resulting object is a once-punctured hyperbolic polygon with a piecewise geodesic boundary of finite length. We refer to such an object as a \textit{polygonal cusp}. For example, Figure~\ref{fig: Cusp}$(ii)$ illustrates such a construction when $m=5$. 
     \begin{figure}[htbp]
    \centering
    \begin{tikzpicture}
        \draw (0,0) circle [radius=3 cm];
        \fill[gray!20!white] 
        ({cos(0)},{sin(0)}) to [bend left] ({cos(60)},{sin(60)}) to [bend right=10] ({3*cos(30)},{3*sin(30)}) to [bend right=10] ({cos(0)},{sin(0)})--cycle;
        \draw ({cos(0)},{sin(0)}) to [bend left] ({cos(60)},{sin(60)});
        \draw[bend right=10,->-=.4] ({cos(60)},{sin(60)}) to  ({3*cos(30)},{3*sin(30)});
        \draw [bend right=10,-<-=.6] ({3*cos(30)},{3*sin(30)}) to ({cos(0)},{sin(0)});

        \fill[gray!20!white] 
        ({cos(0+60)},{sin(0+60)}) to [bend left] ({cos(60+60)},{sin(60+60)}) to [bend right=10] ({3*cos(30+60)},{3*sin(30+60)}) to [bend right=10] ({cos(0+60)},{sin(0+60)})--cycle;

        \draw ({cos(0+60)},{sin(0+60)}) to [bend left] ({cos(60+60)},{sin(60+60)});
        \draw [bend right=10,->-=.4] ({cos(60+60)},{sin(60+60)}) to ({3*cos(30+60)},{3*sin(30+60)});
        \draw [bend right=10,-<-=.6] ({3*cos(30+60)},{3*sin(30+60)}) to [bend right=10] ({cos(0+60)},{sin(0+60)});

        \fill[gray!20!white] 
        ({cos(0+120)},{sin(0+120)}) to [bend left] ({cos(60+120)},{sin(60+120)}) to [bend right=10] ({3*cos(30+120)},{3*sin(30+120)}) to [bend right=10] ({cos(0+120)},{sin(0+120)})--cycle;
        \draw ({cos(0+120)},{sin(0+120)}) to [bend left] ({cos(60+120)},{sin(60+120)});
        \draw [bend right=10,->-=.4] ({cos(60+120)},{sin(60+120)}) to ({3*cos(30+120)},{3*sin(30+120)});
        \draw [bend right=10,-<-=.6] ({3*cos(30+120)},{3*sin(30+120)}) to ({cos(0+120)},{sin(0+120)});

        \fill[gray!20!white] 
        ({cos(0+180)},{sin(0+180)}) to [bend left] ({cos(60+180)},{sin(60+180)}) to [bend right=10] ({3*cos(30+180)},{3*sin(30+180)}) to [bend right=10] ({cos(0+180)},{sin(0+180)})--cycle;
        \draw ({cos(0+180)},{sin(0+180)}) to [bend left] ({cos(60+180)},{sin(60+180)});
        \draw [bend right=10,->-=.4] ({cos(60+180)},{sin(60+180)}) to ({3*cos(30+180)},{3*sin(30+180)});
        \draw [bend right=10,-<-=.6] ({3*cos(30+180)},{3*sin(30+180)}) to ({cos(0+180)},{sin(0+180)});

        \fill[gray!20!white] 
        ({cos(0+240)},{sin(0+240)}) to [bend left] ({cos(60+240)},{sin(60+240)}) to [bend right=10] ({3*cos(30+240)},{3*sin(30+240)}) to [bend right=10] ({cos(0+240)},{sin(0+240)})--cycle;
        \draw ({cos(0+240)},{sin(0+240)}) to [bend left] ({cos(60+240)},{sin(60+240)}) to [bend right=10] ({3*cos(30+240)},{3*sin(30+240)}) to [bend right=10] ({cos(0+240)},{sin(0+240)})--cycle;
        \draw [bend right=10,->-=.4] ({cos(60+240)},{sin(60+240)}) to ({3*cos(30+240)},{3*sin(30+240)});
        \draw [bend right=10,-<-=.6] ({3*cos(30+240)},{3*sin(30+240)}) to ({cos(0+240)},{sin(0+240)});

        \fill[gray!20!white] 
        ({cos(0+300)},{sin(0+300)}) to [bend left] ({cos(60+300)},{sin(60+300)}) to [bend right=10] ({3*cos(30+300)},{3*sin(30+300)}) to [bend right=10] ({cos(0+300)},{sin(0+300)})--cycle;
        \draw ({cos(0+300)},{sin(0+300)}) to [bend left] ({cos(60+300)},{sin(60+300)});
        \draw [bend right=10,->-=.4] ({cos(60+300)},{sin(60+300)}) to  ({3*cos(30+300)},{3*sin(30+300)});
        \draw [bend right=10,-<-=.6] ({3*cos(30+300)},{3*sin(30+300)}) to ({cos(0+300)},{sin(0+300)});

        \draw ({2*cos(20)},{2*sin(20)}) node {$f$} ({2*cos(40)},{2*sin(40)}) node {$a$} ({2*cos(80)},{2*sin(80)}) node {$a$} ({2*cos(100)},{2*sin(100)}) node {$b$} ({2*cos(140)},{2*sin(140)}) node {$b$} ({2*cos(160)},{2*sin(160)}) node {$c$} ({2*cos(200)},{2*sin(200)}) node {$c$} ({2*cos(220)},{2*sin(220)}) node {$d$} ({2*cos(260)},{2*sin(260)}) node {$d$} ({2*cos(280)},{2*sin(280)}) node {$e$} ({2*cos(320)},{2*sin(320)}) node {$e$} ({2*cos(342)},{2*sin(342)}) node {$f$};

        \draw ({cos(0+300)},{sin(0+300)}) to [bend left] ({cos(60+300)},{sin(60+300)}) to [bend right=10] ({3*cos(30+300)},{3*sin(30+300)}) to [bend right=10] ({cos(0+300)},{sin(0+300)})--cycle;

        \draw [bend left, dashed] ({6+1.5*cos(0)},{1.5*sin(0)-1.5}) to ({6+1.5*cos(60)},{1.5*sin(60)-1.5}) to ({6+1.5*cos(120)},{1.5*sin(120)-1.5}) to ({6+1.5*cos(180)},{1.5*sin(180)-1.5});
        
        \draw [bend left] ({6+1.5*cos(180)},{1.5*sin(180)-1.5}) to ({6+1.5*cos(240)},{1.5*sin(240)-1.5}) to ({6+1.5*cos(300)},{1.5*sin(300)-1.5}) to ({6+1.5*cos(360)},{1.5*sin(360)-1.5});
        \draw [] ({6+1.5*cos(0)},{1.5*sin(0)-1.5}) to (6,4-1);
        \draw ({6+1.5*cos(180)},{1.5*sin(180)-1.5}) to (6,4-1);
        \draw [thick, dotted] ({6+1.5*cos(120)},{1.5*sin(120)-1.5}) to (6,4-1);
        \draw [thick, dotted] ({6+1.5*cos(60)},{1.5*sin(60)-1.5}) to (6,4-1);
        \draw  ({6+1.5*cos(240)},{1.5*sin(240)-1.5}) to (6,4-1);
        \draw  ({6+1.5*cos(300)},{1.5*sin(300)-1.5}) to (6,4-1);

        \draw (0,-3.5) node {$(i)$} (6,-3.5) node {$(ii)$};
    \end{tikzpicture}
    \caption{Polygonal Cusp}
    \label{fig: Cusp}
\end{figure}

\begin{figure}[htbp]
    \centering
    \begin{tikzpicture}
        \draw (0,0) circle [radius=2 cm];
        \draw (.2,-.6) to [bend right=20] (-.8,.5);
        \draw (.2,-.6) to [bend left=15] ({2*cos(30)},{2*sin(30)}) to [bend left=15] (-.8,.5);
        \draw (.2,-.25) node {\tiny$\phi$} (-.3,.35) node {\tiny$\theta$} (-.3,-.1) node {\tiny$l$};
        
    \end{tikzpicture}
    \caption{$\Delta_{\theta,l,\phi}$ triangle in Poincar\'e disk $\mathbb{D}$}
    \label{fig:delta triangle}
\end{figure}

The \emph{perimeter} of a polygonal cusp is the sum of the lengths of the geodesic segments of the boundary. If the lengths of the geodesic segments of the boundary are all equal as well as the interior angles at the vertices, then it is called \emph{regular}. 

Let $\Delta_{\theta,l,\phi}$ denote a hyperbolic triangle with exactly one ideal vertex, the length of the side opposite to the ideal vertex $l$ and the interior angles at the other vertices $\theta$ and $\phi$ (see Figure~\ref{fig:delta triangle}). A $p$-sided polygonal cusp can be decomposed into $p$ many triangles $\Delta_{\theta_i,l_i,\phi_i}, i=1, 2, \dots,p,$ where $l_i$'s are the side lengths of the polygonal cusp and $\theta_i, \phi_i$'s are positive real numbers satisfying $\left(\phi_{i-1}+\theta_i\right)$'s are the interior angles of the polygonal cusp (see Figure \ref{fig: 2.3}).

\begin{figure}[htbp]
    \centering
    \begin{tikzpicture}
        \draw (0,0) to [bend left] (2,-.3) to [bend left] (4,0) to [bend left] (5.5,1);
        \draw[dashed] (5.5,1) to [bend left] (6,1.5);
        \draw [dashed] (-1,1) to [bend left] (0,0);
        \draw (0,0) to [bend right=10] (3,5);
        \draw (2,-.3) to [bend left=5] (3,5);
        \draw (4,0) to [bend left=15] (3,5);
        \draw (5.5,1) to [bend left=15] (3,5);
        \draw (1,-.2) node{$l_{i-1}$} (3,-.2) node{$l_{i}$} (4.9,0.4) node{$l_{i+1}$};

        \draw (.7,.3) node {\tiny$\theta_{i-1}$} (1.7,.2) node {\tiny$\phi_{i-1}$} (2.3,.2) node {\tiny$\theta_{i}$} (3.5,.4) node {\tiny$\phi_{i}$} (4.05,.6) node {\tiny$\theta_{i+1}$};


    \end{tikzpicture}
    \caption{Decomposition of a polygonal cusp}
    \label{fig: 2.3}
\end{figure}

\subsection*{Technical lemmas}
Now, we develop two lemmas below, which are essential for the proof of Theorem~\ref{lem: isom}. 

Let $\Sigma_{p, c}$ denote the set of all hyperbolic triangles with perimeter $p$ and a side of constant length $c$, where $p, c\in \mathbb{R}_{>0}, c< p$. The following lemma deals with the maximum possible area of a triangle in $\Sigma_{p, c}$.
\begin{lemma}\label{lem: triangle}
Among the triangles in $\Sigma_{p, c}$, a triangle $\triangle ABC$, with the side $AB$ of fixed length $c$, has the maximum area if and only if the sides $AC$ and $BC$ have equal length.
\end{lemma}

    \begin{figure}[htbp]
        \centering
        \begin{tikzpicture}
            \draw (0,0) to[bend left=15] (4,0) to [bend left=15] (3,2.5) to [bend left=15] (0,0);
            \draw (-.2,-.1) node {$A$} (4.2,-.1) node {$B$}
            (2.8,2.6,-.1) node {$C$};
            \draw (2,0) node{$c$} (1.8,1.4) node{$x$} (3.75,1.4) node{$l-x$};

            \draw (1.3,.45) node {\tiny $\alpha(x)$} (3.3,.4) node {\tiny $\beta(x)$} (2.8,1.7) node {\tiny $\gamma(x)$};
        \end{tikzpicture}
        \caption{}
        \label{fig: 2.2}
    \end{figure}
\begin{proof}
    The perimeter $p$ and the length $c$ of the side $AB$ of a triangle $\Delta ABC\in \Sigma_{p, c}$ are fixed. Therefore, the total length $l$ of the sides $AC$ and $BC$ is fixed. Let us denote the length of the side $AC$ by $x$. Then the length of $BC$ is $l-x$ (see Figure~\ref{fig: 2.2}). Now, applying triangle inequality, we have
    \begin{equation}
        \frac{l-c}{2}<x<\frac{l+c}{2}.
    \end{equation}

  Let $\alpha(x),\beta(x)$ and $\gamma(x)$ denote the interior angles of the triangle $\Delta ABC$ at the vertices $A, B$, and $C$, respectively. Using hyperbolic trigonometric formula (Theorem 2.2.1 \cite{Buser}) to $\triangle ABC$, we get
    \begin{align*}
        \alpha(x)=&\cos^{-1}\left( \frac{\cosh x \cosh c-\cosh(l-x)}{\sinh x \sinh c} \right),\\ \beta(x)=&\cos^{-1}\left( \frac{\cosh(l-x) \cosh c-\cosh x}{\sinh (l-x) \sinh c} \right), \text{ and}\\
        \gamma(x)=&\cos^{-1}\left( \frac{\cosh x \cosh (l-x)-\cosh c}{\sinh x \sinh (l-x)} \right).
\end{align*}
We define $F(x)=\alpha(x)+\beta(x)+\gamma(x), \text{ for } \frac{l-c}{2}<x<\frac{l+c}{2}$. Gauss--Bonnet theorem implies that the area of the triangle $\Delta ABC$ is maximum when $F(x)$ is minimum. Therefore, to complete the proof, we find the point of minima of the function $F$. We have

    \begin{align*}
        \frac{d\alpha}{dx}=&-\frac{1}{\sinh{x}} \sqrt{\frac{\cosh l- \cosh c}{\cosh{c}- \cosh(l-2x)}},\\
         \frac{d\beta}{dx}=&\frac{1}{\sinh(l-x)}\sqrt{\frac{\cosh l - \cosh c}{\cosh c - \cosh(l- 2x)}}, \text{ and }\\
         \frac{d\gamma}{dx}=&\frac{\sinh(l-2x)}{\sinh(l-x) \sinh x}\sqrt{\frac{\cosh l - \cosh c}{\cosh c - \cosh (l-2x)}}.
    \end{align*}

    Now $\frac{dF}{dx}=\left( -\frac{1}{\sinh{x}}+ \frac{1}{\sinh(l-x)}+ \frac{\sinh(l-2x)}{\sinh(l-x) \sinh x}\right) \sqrt{\frac{\cosh l - \cosh c}{\cosh c - \cosh (l-2x)}}.$ The critical points are then determined by 
    $$-\sinh(l-x)+ \sinh x + \sinh (l-2x) =0,$$
    which has a unique zero for $0<x<l$ given by $x=l/2$ (as can be checked by expanding $\sinh$ and solving for $x$). This implies that the function $F(x)$ has a unique minima at $\frac{l}{2}$ and hence the triangle has maximum area only when $AC=BC=\frac{l}{2}$.
\end{proof}

\begin{remark}\label{rem: 2.5}
    The function $F\left(\frac{l}{2}\right)=2\cos^{-1}\left( \coth\frac{l}{2}\tanh\frac{c}{2}\right)+\cos^{-1}\left(\coth^2\frac{l}{2}- \operatorname{cosech}^2\frac{l}{2}\cosh c\right)$ is decreasing with $l$ and hence the area of the isosceles triangle as obtained in Lemma \ref{lem: triangle} is increasing with $l$.
\end{remark}
In the next lemma, we prove the dual statement of Lemma \ref{lem: triangle}. Let $\Sigma^{a, c}$ be the set of all hyperbolic triangles of constant area $a\in \mathbb{R}_{>0}$ and a side of constant length $c$. The following lemma studies the minimal possible perimeter of a triangle in $\Sigma^{a, c}$.

\begin{lemma}\label{dual triangle}
    Among the triangles in $\Sigma^{a, c}$, a triangle $\triangle ABC$, with the side $AB$ of constant length $c$, has the least perimeter only when $AC=BC$.
\end{lemma}
\begin{proof}
    The proof is by contradiction. Suppose $\triangle \in \Sigma^{a,c}$ is a minimum of the perimeter function, for $AC=l_1$ and $BC=l_2$, where $l_1\neq l_2$. Note that such a minima exists because, endowed with the natural topology on $\Sigma^{a,c}$, the perimeter function is proper. Now consider the triangle $\triangle ABC'$ where $AC'=BC'=\frac{l_1+l_2}{2}=\frac{l}{2}$ (say). It follows from Lemma \ref{lem: triangle} that $\mathrm{Area }(\triangle ABC')>\mathrm{Area }(\triangle ABC)$. By Remark \ref{rem: 2.5}, there exists $l'<l$ such that $\mathrm{Area }(\triangle ABC)=\mathrm{Area }(\triangle ABC'')$, where $AC''=BC''=\frac{l'}{2}$. This contradicts the minimality of the perimeter.
\end{proof}

Next we prove the proposition below, which will be used in the proof of Theorem \ref{lem: isom}.

\begin{proposition}\label{prop: 2.6}
    Let $0<a<b\in \mathbb{R}$. For $x>0$, $\Delta_x$ denotes the hyperbolic triangle with vertices at $ai, bi$ and $x$. Define $\mathcal{A}: \mathbb{R}_{>0}\to \mathbb{R}$ by $\mathcal{A}(x)=\mathrm{area}(\Delta_x)$. Then we have the following:
    \begin{enumerate}
        \item The function $\mathcal{A}$ admits its global maximum value at $x=\sqrt{ab}$.
        \item The interior angles of $\Delta_x$ at the vertices $ai$ and $bi$ are the same when $x=\sqrt{ab}.$
    \end{enumerate}
\end{proposition}

\begin{proof}
      Without loss of generality, we assume that $ab=1$. Indeed, if $ab\neq1$, then after applying the isometry $z \mapsto \lambda z$ for an appropriate choice of $\lambda$, we may reduce to the case $ab=1$.
      
    The center $z_0$ and the radius $r$ of the circle passing through $x$ and $ai$ and perpendicular to the real axis are given by
    $$z_0=\frac{x^2-a^2}{2x}, \ \ \ \ \   r=\frac{x^2+a^2}{2x}.$$
    Similarly, the center $z_0'$ and the radius $r'$ of the circle passing through $x$ and $bi$ and perpendicular to the real axis are given by
    $$z_0'
    =\frac{x^2-b^2}{2x}, \ \ \ \ \   r'=\frac{x^2+b^2}{2x}.$$
    Let the interior angles at the vertices $ai$ and $bi$ of the triangle be $\alpha_x$ and $\beta_x$ respectively. By the Gauss--Bonnet theorem, we have
    $$\mathcal{A}(x)=\pi - (\alpha_x+\beta_x).$$
    To find the maximum value of the function, it suffices to find the minimum value of the function $f(x)=\alpha_x+\beta_x$. Now, the following cases may arise.
    
    \noindent \textbf{Case 1:} $x>b$. Then both $z_0$ and $z_0'$ are positive real numbers (see Figure \ref{fig: 0.3}). Also, we have
    $$f(x) = \alpha_x+\beta_x=\frac{\pi}{2}+\sin^{-1}\left( \frac{2ax}{x^2+a^2} \right) + \cos^{-1}\left( \frac{2bx}{x^2+b^2} \right), \  x>b.$$
    For $x>b,\  f(x)$ is differentiable and 
    $$\frac{df}{dx}=\frac{2(x^2-1)(b-a)}{(a^2+x^2)(b^2+x^2)}>0.$$
    Hence, $f$ is strictly increasing in the interval $(b,\infty$).
    \begin{figure}[htbp]
        \centering
        \begin{tikzpicture}
            \draw (-2,0)--(4,0);
            \draw (3,0) arc[start angle=0, end angle=180, radius=1.5cm];
            \draw (3,0) arc[start angle=0, end angle=180, radius=2cm];
            \draw (.5,0)--(.5,3);
            \draw (1.5,0) node {\tiny$\bullet$} (1,0) node {\tiny$\bullet$};
            \draw (3,-.2) node {$x$} (1.5,-.3) node {$z_0$} (1,-.3) node {$z_0'$} (.5,-.3) node {$0$};
            \draw [] (1.5,0)--(.5,{sqrt(1.25)});
            \draw [] (1,0)--(.5,{sqrt(3.75)});
            \draw (.2,{sqrt(1.25)+.1}) node {$ai$} (.2,{sqrt(3.75)+.2}) node {$bi$};

        \end{tikzpicture}
        \caption{}
        \label{fig: 0.3}
    \end{figure}

    \noindent \textbf{Case 2:} $0<x<a$. Then $\frac{1}{a}=b<\frac{1}{x}<\infty$. As $\frac{1}{x}$ is strictly decreasing with $x$, by Case 1 it follows that $f$ is strictly decreasing in $0<x<a$.

    \noindent \textbf{Case 3:} $a \leq x \leq b$. It is straightforward to check $f(a)=f(b)=\frac{\pi}{2}+\sin^{-1}\left( \frac{2}{a^2+b^2} \right)$. Now, we check for $a < x < b$. Then $z_0$ is positive and $z_0'$ is negative (see Figure \ref{fig:0.4}) and 
    $$f(x) = \alpha_x+\beta_x=\sin^{-1}\left( \frac{2ax}{x^2+a^2} \right) + \sin^{-1}\left( \frac{2bx}{x^2+b^2} \right), \  a<x<b.$$
     In $a<x<b,\  f(x)$ is differentiable and 
    $$\frac{df}{dx}=\frac{2(x^2-1)(b-a)}{(a^2+x^2)(b^2+x^2)}.$$
    For $x< 1$, $f'(x)<0$ and for $x> 1$, $f'(x)>0$. Therefore, the function $f$ has a minima at $1$ with
    $$f(1)=2\sin^{-1}\left( \frac{2}{a+b} \right).$$
    As $f(1)<f(a)$, it follows that $x=1$ is the only minima of $f$. This completes the proof of $(1)$. The proof of (2) follows from the following
    $$\alpha_1=\sin^{-1}\left( \frac{2}{a+b} \right)=\beta_{1}.$$
    This completes the proof.
    \begin{figure} [htbp]
        \centering
        \begin{tikzpicture}
            \draw (-4.1,0)--(2.5,0);
            \draw (1.5,0) arc [start angle=0, end angle=180, radius=1];
            \draw (1.5,0) arc [start angle=0, end angle=180, radius=2.3];
            \draw (0,0)--(0,3);
            \draw (.5,0)--(0,{sqrt(.75)});
            \draw (-.8,0)--(0,{sqrt(4.65)});
            \draw (.5,-.3) node {$z_0$} (-.8,-.3) node {$z_0'$} (0,-.3) node {$0$} (1.5,-.3) node {$x$};
            \draw (.5,0) node {\tiny$\bullet$} (-.8,0) node {\tiny$\bullet$};
            \draw (.25,{sqrt(4.65)+.25}) node {$bi$} (.25,{sqrt(.75)+.35}) node {$ai$};
        \end{tikzpicture}
        \caption{}
        \label{fig:0.4}
    \end{figure}
\end{proof}
\begin{corollary} \label{cor 0.10}
    If the length of the finite side of the triangle $\Delta_x$, where $x=\sqrt{ab}$ is $l$, then
    $$l=\log\frac{b}{a}\implies \frac{b}{a}=e^l.$$
    Therefore, each of the interior angles (non-ideal) of the triangle with maximum area is 
    $$\alpha_l= \sin^{-1}\left( \frac{2\sqrt {ab}}{a+b} \right)= \sin^{-1}\left( \frac{2e^{\frac{l}{2}}}{1+e^l} \right) = \sin^{-1}\left(\operatorname{sech}\left(\frac{l}{2}\right)\right).$$
    Note that $\alpha_l \longrightarrow 0$ as $l\longrightarrow \infty$ and $\alpha_l \longrightarrow \frac{\pi}{2}$ as $l\longrightarrow 0$ and $\alpha_l$ is strictly decreasing for $l>0$. Therefore, there exists a unique $l_0>0$ such that $\alpha_{l_0}=\frac{\pi}{4}$.
\end{corollary}

\begin{proof}[Proof of Theorem \ref{lem: isom}]
    A polygonal cusp can be decomposed into $\Delta_{\theta,l,\phi}$ triangles. Take two such adjacent triangles which have a common edge, for example, $\triangle ABO$ and $\triangle BCO$ in Figure \ref{fig: 2}. By Lemma \ref{lem: triangle}, the area of the quadrilateral ABCO is the maximum only when $AB=BC$. That the angles of the polygonal cusp are also equal follows from Proposition \ref{prop: 2.6}. This completes the proof.
 \begin{figure}[htbp]
        \centering
        \begin{tikzpicture}
            \draw (2,0) to [bend left=25] (0,4);
            \draw (-2,0) to [bend right=25] (0,4);
            \draw (-2,0) to [bend left=50] (-.5,-1) to [bend left=50] (2,0);
            \draw  [line width=0.1pt, opacity=0.2] (-2,0) to [bend left=30] (2,0);
            \draw (0,4) to [bend left=5] (-.5,-1);
            \draw (.3,4) node {$O$} (-2.3,0) node{$A$} (-.25,-1) node {$B$} (2.3,0) node{$C$};
        \end{tikzpicture}
        \caption{}
        \label{fig: 2}
    \end{figure}
\end{proof}

Suppose that the numbers $\theta$ and $l$ satisfy  
$$\theta=\sin^{-1}\left(\operatorname{sech} (l/2)\right).$$

The symbol $\Delta_{l,\theta}$ denotes a triangle $\Delta_{\theta,l,\theta}$.

Now, we state and prove the dual statement of Theorem \ref{lem: isom}.
\begin{theorem}\label{dual isom}
Among all $p$-sided polygonal cusps with a fixed hyperbolic area, the one with the least perimeter is the regular one. 
\end{theorem}
    \begin{proof}
        It follows from Lemma \ref{dual triangle} that the side lengths of the polygonal cusp admitting the minimal perimeter are all equal, say $l$. Decompose it into $\Delta_{\theta,l,\phi}$-triangles and let $\Delta_{\theta_1,l,\theta_2}$ be such a triangle. By Proposition \ref{prop: 2.6}, the area of the triangle $\Delta_{l,\theta}$ is greater than that of $\Delta_{\theta_1,l,\theta_2}$ and hence $2\theta<\theta_1+\theta_2$. It is clear that $\mathrm{Area}(\Delta_{\overline{l},\overline{\theta}})=\mathrm{Area}(\Delta_{\theta_1,l,\theta_2})<\mathrm{Area}(\Delta_{l,\theta})$, where $\overline{\theta}=\frac{\theta_1+\theta_2}{2}$. As $\theta< \overline{\theta}$, by Corollary \ref{cor 0.10}, we have $l>\overline{l}$. But this contradicts the minimality of the perimeter.
    \end{proof}

\begin{definition}
    A regular polygonal cusp with each side length $l$ and interior angle $2\theta$ is called a \textit{polygonal cusp of type $\Delta_{l,\theta}$}.
\end{definition}

\section{Length of filling pairs on once punctured surfaces}
The goal of this section is to prove Theorem \ref{thm 1.1}. The idea of the proof is adapted from the work of Aougab--Huang (Theorem 1.3 \cite{Aougab}) and Gaster's proof \cite{Gaster} to a conjecture of Aougab-Huang on the length of filling pairs on closed surfaces. Before going to the proof, we recall the definition of spreads in a graph embedded on an oriented surface and develop a key lemma that is used in the proof of Theorem~\ref{thm 1.1}. 

Let $G$ be a graph embedded on an oriented surface $S$. The orientation of the surface gives a cyclic ordering to the edges incident at $p$, for every vertex $p$ of the graph $G$. A sub-graph $H\subset G$ is called a \emph{spread} if for every vertex $p$ in $H$ and edges $e,e'$ of $H$ incident at $p$, the edges $e$ and $e'$ are not consecutive in the cyclic order at $p$. A sub-graph $H$ of $G$ is called spanning if they have the same vertex set and a forest if $H$ has no cycle. A forest $H$ is called a tree if it is connected.  

Given $(\alpha, \beta)$ a filling pair on a surface $S_{g, n}$, we may regard this as a graph $\Gamma(\alpha, \beta)$ on the surface, where the intersection points are the vertices and the sub-arcs of $\alpha$ and $\beta$ between the vertices are the edges. 

The main ingredient in Gaster's proof \cite{Gaster} of the Aougab--Huang conjecture (Section 4, \cite{Aougab}) on the length of filling pairs on closed surfaces is the existence of a spread spanning tree and a spread spanning forest in the dual graph of $\Gamma(\alpha,\beta)$ associated with a filling pair $(\alpha,\beta)$.
 Our arguments for proving Theorem \ref{thm 1.1} follow Gaster’s overall strategy. However, a direct adaptation to the punctured setting introduces additional complications, since some vertices of the dual graph of $\Gamma(\alpha,\beta)$ correspond to punctured disks, unlike in the closed-surface case. To address this issue, we adopt the following convention: we regard the punctures as marked points on the surface and define these marked points to be the vertices of the dual graph corresponding to the punctured disks. The spread spanning tree or forest is used only to determine the order in which the decomposed polygons and polygonal cusps are attached to obtain either a single polygonal cusp or a polygon together with a polygonal cusp and therefore the convention serves our purpose. Next, we obtain the following lemma.

\begin{lemma}\label{lem 3.1}
    Let $(\alpha, \beta)$ be a filling pair on a surface $S_{g,n}$. If $\alpha$ is non-separating, then the dual graph of $\alpha\cup \beta$ admits a spread spanning tree, and if $\alpha$ is separating, the dual graph admits a spread spanning forest with two components.
\end{lemma}

\begin{proof}
After treating the punctures as marked points, the problem becomes identical to the corresponding one in the closed-surface case. Therefore, the existence of a spread spanning tree or spread spanning forest follows exactly as in Gaster's proof, and we omit the details.
\end{proof}
Now, we are ready to prove Theorem \ref{thm 1.1}.

\begin{proof}[Proof of Theorem \ref{thm 1.1}]
    For a minimally intersecting filling pair $(\alpha, \beta)$ on a hyperbolic surface $X\in \mathcal{T}_{g,1}$, the complement  $X\setminus (\alpha\cup\beta)$ is a $(8g-4)$-sided polygonal cusp. The area of the polygonal cusp is the same as the area of the surface $X$, which is $-2\pi \chi(X)=2\pi(2g-1)$. On the other hand, if $\Sigma$ denotes the sum of the interior angles of the polygonal cusp, then by Gauss--Bonnet theorem, we have
    \begin{align*}
    (8g-4)\pi-\Sigma & =2\pi(2g-1)) \\ \implies \Sigma &=2\pi(2g-1).
    \end{align*}
    
    By Theorem \ref{dual isom}, a regular polygonal cusp gives the least perimeter for a given fixed area, and hence each interior angle of the polygonal cusp is $\frac{2\pi(2g-1)}{8g-4}=\frac{\pi}{2}$. As discussed in the proof of Theorem \ref{dual isom}, the polygonal cusp with minimal perimeter can be decomposed into triangles, each of type $\Delta_{l,\frac{\pi}{4}}$, where $l$ is given by 
    \begin{align*}
\frac{\pi}{4} &=\sin^{-1} \left( \frac{2e^{l/2}}{1+e^l} \right)\\ \implies l &=2\ln{\left(\sqrt{2}+1\right)}.
\end{align*}
This proves the theorem for minimally intersecting filling pairs, and it also follows that the inequality is sharp.

Now suppose that $(\alpha, \beta)$ is any filling pair on $X$. We consider the cases below:

\noindent \textbf{Case 1:} At least one of the curves, say $\alpha$ is non-separating. Let $$X\setminus(\alpha\cup \beta)=\bigsqcup_{k=1}^rP_k,$$ where for every $k\in \{1, 2, \dots, r\}$, $P_k$ is either a hyperbolic polygon or a polygonal cusp. The filling pair $(\alpha, \beta)$ gives a cellular decomposition of the surface $X$ with number of 0-cells $i(\alpha,\beta)$, number of 1-cells $2i(\alpha, \beta)$. By Euler characteristic arguments, we have
    $$i(\alpha,\beta)=2g-2+r.$$
    Valency condition of the graph $\alpha\cup\beta$ on $S_g$ implies that
\begin{align*}
    \Sigma_{k=1}^r n\left(P_k\right) = 4i(\alpha,\beta) = 8g-8+4r,
\end{align*}
where $n\left(P_k\right)$ denotes the number of vertices in $P_k$. By Lemma \ref{lem 3.1}, the dual of the graph $\alpha\cup \beta$ admits a spread spanning tree $T$. We attach the polygons $P_k$ using the tree $T$ to obtain a polygonal cusp $Q$ as follows. Let $e$ be an edge of $T$ whose endpoints are the polygons $P_i$ and $P_j$. Attach these two polygons along the edge dual to $e$. We perform this operation for all edges of $T$ and obtain a polygonal cusp, $\overline{Q}$. Some vertices of $\overline{Q}$ are the endpoints of some edges of $P_k$'s along which they are identified. We call such a verex as a \emph{new-vertex}. As $T$ is a spread, the angle at each of these new vertices is $\pi$. Therefore, we may forget the new vertices and call the polygonal cusp as $Q$. Now, we calculate the number of vertices of $Q$. For each edge $e$ of $T$, exactly four vertices of $\bigsqcup_{k=1}^rP_k$ project onto some new-vertices. Since $T$ has $r-1$ edges,
    $$n(Q)=\Sigma_{k=1}^rn\left(P_k\right)-4(r-1)=8g-4.$$
    Thus, $Q$ is a $8g-4$ sided polygonal cusp whose area is $2\pi(2g-2)$, so its perimeter is at least $2(8g-4)\ln(\sqrt{2}+1)$. This implies that 
    $$l(\alpha,\beta)=\frac{1}{2}\cdot \Sigma_{k=1}^r\mathrm{Perim}\left(P_k\right)\geq \frac{1}{2}\cdot \mathrm{Perim}(Q)\geq(8g-4)\ln(\sqrt{2}+1).$$
    \noindent \textbf{Case 2:} Both $\alpha$ and $\beta$ are separating. Let $X\setminus\alpha=X_{g_1}\bigsqcup X_{g_2}$, where $X_{g_1}\in \mathcal{T}_{g_1,1}$, $X_{g_2}\in \mathcal{T}_{g_2}$ and $g_1+g_2=g$. Let $X_{g_j}\setminus\beta=\bigsqcup_{k=1}^{r_j}P_k^{(j)},j=1,2$. By Lemma \ref{lem 3.1}, the dual graph of $\alpha\cup \beta$ admits a spread spanning forest $T_1\bigsqcup T_2$, where $T_j$ is a spread spanning tree of the dual of $(\alpha\cup \beta)\cap X_{g_j},j=1,2$. Using the Euler characteristic arguments as in Case 1, we have 
    $$\Sigma_{k=1}^r n\left(P_k^{(j)}\right)=4i(\alpha,\beta)= 8g_j-4+4r_j,\ \ j=1,2.$$
    As in Case 1, we construct a polygonal cusp $Q_1$ and a polygon $Q_2$ using $T_1$ and $T_2$ respectively, with the number of vertices $n(Q_j)=8g_j, j=1,2$. It is clear that 
    $$l(\alpha,\beta)\geq \frac{1}{2}\mathrm{Perim}(Q_1)+\frac{1}{2}\mathrm{Perim}(Q_2).$$
    The remaining part follows from Bezdek's \cite{Bezdek} isoperimetric inequality for hyperbolic polygons and  Lemma \ref{lemma 3.2} below.
\end{proof}
 \begin{lemma}\label{lemma 3.2}
        Let $R_1$ be a $n_1$ sided regular right-angled hyperbolic polygon; $R_2$ and $R$ be respectively $n_2$ and $m$ sided polygonal cusps with $n_1+n_2=m+4$ and both of them are of type $\Delta_{l,\theta}$. Then
        $$\mathrm{Perim}(R_1) +\mathrm{Perim}(R_2)\geq\mathrm{Perim}(R).$$
    \end{lemma}
    \begin{proof}
        Consider two functions $f(x)=2x\cosh^{-1}\left( \sqrt{2}\cos\frac{\pi}{x}\right)$, $g(x)=2x\ln(\sqrt{2}+1)$. Then it is clear that $\mathrm{Perim}(R_1)=f(n_1), \mathrm{Perim}(R_2)=g(n_2)$ and $\mathrm{Perim}(R)=g(m)$. In addition, $m,n_2\geq3$ and the Gauss-Bonnet theorem imply that $n_1\geq5$. Observe that, as $g$ is linear, $g(n_2) - g(m) = -g(m - n_2) = -g(n_1 - 4)$. Hence in order to find that $f(n_1)+g(n_2)\geq g(m)$, it suffices to show that $f(n_1)-g(n_1-4)\geq0$.
        Define a function
        \begin{align*}
            h(x)\coloneqq &f(x)-g(x-4)\\
            =&2x\left(\cosh^{-1}\left( \sqrt{2}\cos\frac{\pi}{x}\right)-\ln(\sqrt{2}+1)\right)+8\ln(\sqrt{2}+1)\\
            =&2x\left(\ln\left( \sqrt{2}\cos\frac{\pi}{x}+ \sqrt{\cos\frac{2\pi}{x}}\right)-\ln(\sqrt{2}+1)\right)+8\ln(\sqrt{2}+1).
        \end{align*}
        \begin{align*}
            h'(x)=&2\left(\ln\left( \sqrt{2}\cos\frac{\pi}{x}+ \sqrt{\cos\frac{2\pi}{x}}\right)-\ln(\sqrt{2}+1)\right)+\frac{2\pi}{x}\frac{\sqrt{2}\sin\frac{\pi}{x}+\frac{\sin\frac{2\pi}{x}}{\sqrt{\cos\frac{2\pi}{x}}}}{\sqrt{2}\cos\frac{\pi}{x}+ \sqrt{\cos\frac{2\pi}{x}}}\\
            h''(x)=&f''(x)=-\frac{2\pi^2\sqrt{2}\cos\frac{\pi}{x}}{x^3\sqrt{\cos^3\frac{2\pi}{x}}}.
        \end{align*}
        Note that $h''(x)<0$ for $x>4$ and hence $h'(x)$ is strictly decreasing for $x>4$. As $h'(x)\longrightarrow0$ as $x\longrightarrow\infty$, it follows that $h'(x)>0$ for $x>4$ and this implies that $h$ is increasing for $x>4$. As $h(4)=0$, it follows that $h(x)>0$ for all $x>4$. This completes the proof. 
    \end{proof}

\begin{remark}[Surfaces with more than one puncture] An interesting direction for future research is to determine the minimum length of filling pairs on hyperbolic surfaces with more than one puncture. The methods used in this paper do not seem to generalize to this setting. Indeed, the proof would require gluing two polygonal cusps, after which Theorem~\ref{lem: isom} is no longer applicable. It would therefore be interesting to develop new techniques to handle this case.
\end{remark}
\section*{Acknowledgments}
The authors would like to thank the anonymous referee for the useful suggestions and comments that helped improve the exposition.

\bibliographystyle{alpha}
\bibliography{bibliography.bib}

@book {Buser,
    AUTHOR = {Buser, Peter},
     TITLE = {Geometry and spectra of compact {R}iemann surfaces},
    SERIES = {Modern Birkh\"auser Classics},
      NOTE = {Reprint of the 1992 edition},
 PUBLISHER = {Birkh\"auser Boston, Ltd., Boston, MA},
      YEAR = {2010},
     PAGES = {xvi+454},
      ISBN = {978-0-8176-4991-3},
   MRCLASS = {58J50 (30F10 32G15 58J53)},
  MRNUMBER = {2742784},
       DOI = {10.1007/978-0-8176-4992-0},
       URL = {https://doi.org/10.1007/978-0-8176-4992-0},
}

@article{thurstonspine,
  title={A spine for Teichmuller space},
  author={Thurston,W.},
  journal={preprint},
  year={1986},
}

@article {Anderson,
    AUTHOR = {Anderson, James W. and Parlier, Hugo and Pettet, Alexandra},
     TITLE = {Small filling sets of curves on a surface},
   JOURNAL = {Topology Appl.},
  FJOURNAL = {Topology and its Applications},
    VOLUME = {158},
      YEAR = {2011},
    NUMBER = {1},
     PAGES = {84--92},
      ISSN = {0166-8641,1879-3207},
   MRCLASS = {57M50},
  MRNUMBER = {2734699},
       DOI = {10.1016/j.topol.2010.10.007},
       URL = {https://doi.org/10.1016/j.topol.2010.10.007},
}

@article {Fanoni,
    AUTHOR = {Fanoni, Federica and Parlier, Hugo},
     TITLE = {Filling sets of curves on punctured surfaces},
   JOURNAL = {New York J. Math.},
  FJOURNAL = {New York Journal of Mathematics},
    VOLUME = {22},
      YEAR = {2016},
     PAGES = {653--666},
      ISSN = {1076-9803},
   MRCLASS = {57M50 (30F10 30F45)},
  MRNUMBER = {3548116},
MRREVIEWER = {Colin\ C.\ Adams},
       URL = {http://nyjm.albany.edu:8000/j/2016/22_653.html},
}

@article {Aougab,
    AUTHOR = {Aougab, Tarik and Huang, Shinnyih},
     TITLE = {Minimally intersecting filling pairs on surfaces},
   JOURNAL = {Algebr. Geom. Topol.},
  FJOURNAL = {Algebraic \& Geometric Topology},
    VOLUME = {15},
      YEAR = {2015},
    NUMBER = {2},
     PAGES = {903--932},
      ISSN = {1472-2747,1472-2739},
   MRCLASS = {57M20 (57M50 57N05)},
  MRNUMBER = {3342680},
MRREVIEWER = {John\ G.\ Ratcliffe},
       DOI = {10.2140/agt.2015.15.903},
       URL = {https://doi.org/10.2140/agt.2015.15.903},
}

@article {Penner,
    AUTHOR = {Penner, Robert C.},
     TITLE = {A construction of pseudo-{A}nosov homeomorphisms},
   JOURNAL = {Trans. Amer. Math. Soc.},
  FJOURNAL = {Transactions of the American Mathematical Society},
    VOLUME = {310},
      YEAR = {1988},
    NUMBER = {1},
     PAGES = {179--197},
      ISSN = {0002-9947,1088-6850},
   MRCLASS = {57N05 (20F34 58F15)},
  MRNUMBER = {930079},
MRREVIEWER = {William\ Dunbar},
       DOI = {10.2307/2001116},
       URL = {https://doi.org/10.2307/2001116},
}

@article {Gaster,
    AUTHOR = {Gaster, Jonah},
     TITLE = {A short proof of a conjecture of {A}ougab-{H}uang},
   JOURNAL = {Geom. Dedicata},
  FJOURNAL = {Geometriae Dedicata},
    VOLUME = {213},
      YEAR = {2021},
     PAGES = {339--343},
      ISSN = {0046-5755,1572-9168},
   MRCLASS = {51M09 (51M25)},
  MRNUMBER = {4278331},
       DOI = {10.1007/s10711-020-00584-w},
       URL = {https://doi.org/10.1007/s10711-020-00584-w},
}

@article {Sanki,
    AUTHOR = {Sanki, Bidyut and Vadnere, Arya},
     TITLE = {A conjecture on the lengths of filling pairs},
   JOURNAL = {Geom. Dedicata},
  FJOURNAL = {Geometriae Dedicata},
    VOLUME = {213},
      YEAR = {2021},
     PAGES = {359--373},
      ISSN = {0046-5755,1572-9168},
   MRCLASS = {57K20 (51M16 51M25)},
  MRNUMBER = {4278333},
MRREVIEWER = {Hengnan\ Hu},
       DOI = {10.1007/s10711-020-00586-8},
       URL = {https://doi.org/10.1007/s10711-020-00586-8},
}

@article {Bezdek,
    AUTHOR = {Bezdek, K.},
     TITLE = {Ein elementarer {B}eweis f\"ur die isoperimetrische
              {U}ngleichung in der {E}uklidischen und hyperbolischen
              {E}bene},
   JOURNAL = {Ann. Univ. Sci. Budapest. E\"otv\"os Sect. Math.},
  FJOURNAL = {Annales Universitatis Scientiarum Budapestinensis de Rolando
              E\"otv\"os Nominatae. Sectio Mathematica},
    VOLUME = {27},
      YEAR = {1984},
     PAGES = {107--112},
      ISSN = {0524-9007},
   MRCLASS = {52A40 (51M25)},
  MRNUMBER = {823098},
MRREVIEWER = {G.\ D.\ Chakerian},
}

@book{farb2011primer,
  title={A primer on mapping class groups (pms-49)},
  author={Farb, Benson and Margalit, Dan},
  year={2011},
  publisher={Princeton university press}
}

\end{document}